\documentclass[11pt]{amsart}
\usepackage{amssymb,epsfig}
\usepackage{amsmath}
\newtheorem{theorem}{Theorem}[section]

\newtheorem{lemma}[theorem]{Lemma}

\newtheorem{remark}[theorem]{Remark}

\usepackage{amsmath,amsfonts,epsfig,latexsym}
\def\square{\hbox{\vrule\vbox{\hrule\phantom{o}\hrule}\vrule}}

{\vspace{2ex}\noindent{\sl Proof: }}{{\null\hfill\square\bigskip}}
\def\O{{\rm O}}
\def\po{{\rm o}}

\def\im{{\rm Im}}
\def\limsup {{\rm limsup}}
\def\Supp {{\rm supp}}
\def\R{\mathbb {R}}
\def\C{\mathbb {C}} 
\def\N{\mathbb {N}}

\newcommand{\limSup}{\mathop{\limsup}}
\begin{document}
\title{Asymptotics of resonances for a schr\"odinger
operator  with matrix values }
\author{ L. Nedelec$\,^\star$}
\address{L. Nedelec,  L.A.G.A., Institut Gali\-l\'ee, Univer\-sit\'e de Paris Nord, av. J.B. Cle\-ment,
F-93430 Villetaneuse, France IUFM de l'academie de Rouen, France}
\email{ nedelec@math.univ-paris13.fr}
\thanks{$\,^\star$  }

\maketitle

\section{Introduction }
We obtain the asymptotics for the counting function of resonances for a matrix valued Schr\"odinger
operator in dimension one. This  theorem  is a generalization of the result of
Zworski \cite{macas}, using ideas from Froese
\cite{froese}, both of which only treat the scalar case. \par

There has been much work in proving upper and lower bounds for the number of
resonances in various situations (for a review of this subject see
\cite{remac}), but only \cite{macas} \cite{macas2},\cite{froese} obtain
results about
asymptotic. Here we generalize these results to matrix potentials which appear in such
physically important models as the Born Oppenheimer approximation, Dirac operator etc. 
Resonances are defined as poles of the continuation of  the scattering matrix.
In physics, the main interest is in resonances close to the real axis, since these
can be measured in experiments. From a mathematical point of view, their
asymptotics are important in order to establish  trace formula.
\par
Compared to \cite{froese}, where the resolvent was used and the potential $V$
was a scalar, the proof here uses the determinant of the scattering
matrix to  characterize
the resonances. We use the Born approximation to
approximate  the scattering matrix. 
 This approximation, based on  WKB solution, is widely used in inverse
 scattering theory and is a good tool here as well.
\par
For Schr\"odinger operators in dimensions greater than one, sharp upper bounds
on the number of resonances are known \cite{mac}, but only weak lower bounds
exist \cite{masa}.
  
\par
The result for the scalar case is
\begin{theorem}\label{mainvieux}
Apart from a set of density zero for the Lebesgue measure, all  resonances of $-\Delta+V$ are contained in arbitrarily small sectors about the real axis. 
Let $n^\pm(r)$ denote the number of resonances of modulus less than $r$ contained in some sector about $\R^\pm$ then
$$n_\pm(r)=C(V)r+\po(r).$$
$$C(V)=\frac 1 \pi \big(\sup_{x,y\in\Supp V}|x-y| \big).$$
\end{theorem}
The quantity $C(V)$ get no clear generalization in the matrix case.\par

We study the resonances of the operator $-\Delta+V$ in dimension one, where $V$ is a matrix with value $V(x)\in M_n(\C)$ for $x\in\R$,  
real  $\overline V =V$  and of compact support and bounded
( so $V\in L^1$).
Define
$$t_{i,j}=\inf_{x\in\Supp V_{ij}} x, \,u_{i,j}= \sup_{x\in\Supp V_{ij}} x,$$

We obtain 
\begin{theorem}\label{main}
Apart from a set of density zero for  Lebesgue measure, all  resonances of $-\Delta+V$ are contained in arbitrarily small sectors about the real axis. 
Let $n^\pm(r)$ denote the number of resonances of modulus less than $r$ contained in some sector about $\R^\pm$ then
 they exist a constant depending only on $V$ such that 
 $$n_\pm(r)=C(V)r+\po(r),$$ 
 $$C(V)\leq  \frac 1 \pi \big(\sum_{i} \sup_{\sigma\in S_n} \sup_{p,p'}u_{i,p}-t_{p',\sigma(i)}\big).$$
In some case the constant $C(V)$ could be determined,

{\bf (H1)} : if $V$ is a triangular matrix  then $$C(V)= \frac 1 \pi \big(\sum_{i} u_{i,i}-t_{i,i}\big).$$
{\bf (H2)} :  if $V$ is a such that for all $i,j$ $V_{i,j}\geq 0$, 
$$\sup_{\sigma\in S_n} \sum_{i} \sup_{p}t_{i,p}-u_{p,\sigma(i)}$$ is attained for only one $\sigma_0\in S_n$ 
and
$$\sup_{\sigma\in S_n}\sum_{i}  \sup_{p\neq p'}\sum_{i} t_{i,p}-u_{p',\sigma(i)}<\sum_{i} \sup_{p}t_{i,p}-u_{p,\sigma_0(i)}$$
and 
$$\sup_{p}\sum_{i} t_{i,p}-u_{p,\sigma_0(i)}$$ is attained for only one $p$, denote by $p(i)$, 
 then 
$$C(V)= -\sup_{\sigma\in S_n} \sum_{i} \sup_{p}\sum_{i} t_{i,p}-u_{p,\sigma(i)}.$$
\end{theorem}

The plan of the article is as follow. 
The two first sections recall known results which we will use.  The first 
 estimates the number of zero of an analytic function already used in
 \cite{froese}, and it allowed him to give asymptotic when $V$ is scalar  but
 does not have compact support. 
The second is essentially the Paley Wiener Theorem. 
Section \ref{wkb} concerns the approximation  of the transmission matrix (Born Approximation).
Section  \ref{trans} is the definition of the transmission matrix.
Section \ref{syme} just  recall some simplifications due to the
symmetry of the potential. 
The section \ref{scatte} contain the proof using all the tools of the previous
sections.    

\section{Zero of entire functions}
Let denote by $S(k)$ the scattering matrix.
The zeroes of $\det S(k)$ are the complex conjugates of the resonances of our
problem so it suffice to count the number of zeroes of  $\det S(k)$ in the
upper half plane. The function   $\det S(k)$ is analytic in the upper half
plane so we going to use Theorem \ref{levin}. 

A function is said to be of exponential type in the upper half plane if 
$$\limSup_{r\to\infty} \sup_{|z|=r, \Im z\geq 0}\frac { \ln|F(z)|}{ r}< \infty .$$  We call is type the number $\limSup_{r\to\infty} \sup_{|z|=r, \Im z\geq 0}\frac { \ln|F(z)|}{ r}$. Denote by $n(r,\theta_1,\theta_2)$ the number of zeroes of $F$ in the sector 
$$\{z;|z|\leq r, \theta_1\leq \im z\leq \theta_2\}.$$
The following result appears in  \cite{levin} p 243 and p 251.

\begin{theorem} \label{levin}
Let $F$ be an holomorphic function on the half plane $\{z, \arg z\geq 0\}$ of
exponential type which satisfies
$$\int_{-\infty}^{\infty}\frac {\ln_+|F(x)| }{1+x^2}dx< \infty.$$
Then the function $\tau(\theta)$ defined by  $\tau(\theta)=\limSup_{r\to\infty} \ln|F(re^{i\theta})|/r$
 satisfies
$$\tau(\theta)=\tau |\sin \theta| \mbox{ for } \theta\in]0,\pi[,$$ where
$\tau$ is constant.
For  $\theta_1,\theta_2\in ]0,\pi[,$ 
$$\lim_{r\to\infty} \frac {n(r,\theta_1,\theta_2)}{r}=\frac 1 {2\pi}\big[\tau'(\theta_1)-\tau'(\theta_2)+\int_{\theta_2}^{\theta_1} \tau(\theta)d\theta\big]=0.$$
For $\theta>0$ 
$$\lim_{r\to\infty} \frac {n(r,0,\theta)}{r}=\frac {\tau} {\pi}.$$
\end{theorem}

\section{An variant of Paley-Wiener}
We quote another lemma from \cite{froese}, that we need,  which is a variant of the Paley-Wiener Theorem.
\begin{lemma}\label{paley}
Assume $W\in L^\infty$, and $\Supp\, W\subset[ -1,1]$ but $\Supp \, W$  contain in no smaller interval. Let $f(k,x)$ be holomorphic in $k$ in the lower half plane, and satisfies uniformly in $x$ 
$$    |f(k,x)|\leq C/|k|.$$ 
Then $\int_{-1}^1 e^{\pm ikx} W(x)(1+f(k,x))  $ is of type at least $1$ in the lower half plan.
\end{lemma}

\section{Approximation using WKB solutions}\label{wkb}
To apply the theorem \ref{levin} to the function $\det S(k)$, we are going to give an approximation of the scattering matrix by computing the transmission matrix to compute $\tau$, and prove that $\det S(k)$ satisfies the hypotheses of Theorem \ref{levin}.
The following considerations will give asymptotics in $k$ of  solutions of the equation 
\begin{equation}\label{eq:gen0}(-\Delta+V-k^2)u=0,\end{equation}
 or 
 equivalently, for the solutions of the system,  
\begin{equation}\label{eq:gen1}
 \left(\begin{array}{cc}
{\frac 1 i \frac d {dx}}&0\\
0&{\frac 1 i \frac d {dx}}
\end{array}\right)  v
=\left(\begin{array}{cc} -k+\frac V {2k}&-\frac V {2k} \\   \frac V {2k}&k-\frac V {2k}\end{array}\right)
 v.
\end{equation}
where $v_1= -iku+\frac d {dx} u$ and $v_2= iku+\frac d {dx} u,$ which is
a convenient reformulation for studying the behavior at infinite of the solution.

Denote the standard basis of $\R^{2n}$ by   $\{e_j\}_{j\in\{1,\cdots,2n\}}$,
and set  $v_j^{-}=e^{- ikx} e_j$ and $v_j^{+}=e^{ ikx} e_{j+n}$ for $j\in\{1,\cdots,n\}$, these form a  base of  the solution of the system
\begin{equation}\label{eq:gen2}
 \left(\begin{array}{cc}
{\frac 1 i \frac d {dx}} &0\\
0&{\frac 1 i \frac d {dx}}
\end{array}\right)  \,v
=\left(\begin{array}{cc} -k &0 \\ 0  &k \end{array}\right) v.
\end{equation}
The function $v=\sum_{j=1}^n \beta^-_j v^-_j+ \beta^+_j v^+_j$ is a solution of the system
(\ref{eq:gen1}) if and only if  
\begin{equation}\label{eq:gen3}
 \left(\begin{array}{cc}
{1 \over i}{ d\over dx}&0\\
0&{1\over i }{d\over dx}
\end{array}\right)  \left(\begin{array}{c} \beta^- \\ \beta^+\end{array}\right)  \,
=\left(\begin{array}{cc} \frac V {2k} &-e^{2ik}\frac V {2k} \\
 e^{-2ik}\frac V {2k}& -\frac V {2k}  \end{array}\right)
\left(\begin{array}{c} \beta^- \\ \beta^+\end{array}\right).
\end{equation}
Writing
$$\gamma^{+}_j=\beta^{+}_j e^{2ikx}\quad \gamma^{-}_j=\beta^{-}_j ,$$
\begin{equation}\label{eq:gen4}
 \left(\begin{array}{cc}
{1\over i}{ d\over dx}&0\\
0&{1\over i}{ d\over dx}-2k
\end{array}\right)  \left(\begin{array}{c} \gamma^- \\ \gamma^+\end{array}\right)  \,
=\left(\begin{array}{cc} \frac V {2k} &-\frac V {2k} \\
\frac V {2k}& -\frac V {2k}\end{array}\right)
\left(\begin{array}{c} \gamma^- \\ \gamma^+\end{array}\right).
\end{equation}

We also write this  in expanded form as 

\begin{equation}\label{eq:gen5}
\begin{array}{c}
{d\over idx}  \gamma^-_j =
 \sum_k \frac {V_{jk}} {2k}   \gamma^-_k
+\sum_k -\frac {V_{jk}} {2k}  \gamma^+_k,\\
{d\over idx}  \gamma^+_j-2k \gamma^+_j =
\sum_k \frac {V_{jk}} {2k}   \gamma^-_k
 +\sum_k -\frac {V_{jk}} {2k} \gamma^+_k.
\end{array}
\end{equation}

Let $x_0\in\overline\R$ such that $x_0>\sup\{x,x\in \Supp V\}$.
This system (\ref{eq:gen5})  can be solved by  iteration
as follows.  Let $\gamma^-_j=\sum_{n=0}^\infty  \gamma^-_{j,n},\,$
 $\gamma^+_j=\sum_{n=0}^\infty  \gamma^+_{j,n}$ , where we  define an initial condition

\begin{equation}\label{eq:gen7}
\begin{array}{c}\gamma^+_{j,0}(x)=0 \quad  \gamma^-_{j,0}(x)=c_j\\
\mbox{ or } \gamma^+_{j,0}(x)=d_je^{2ikx} \quad  \gamma^-_{j,0}(x_0)=0\\
\end{array}
\end{equation}

and an iteration procedure
$$ \gamma^-_{j,n+1}(x)=J^-_j(\gamma^-_{n}-\gamma^+_{n})(x),$$
$$ \gamma^+_{j,n+1}(x)=J^+_j(\gamma^-_{n}-\gamma^+_{n})(x),$$
where $J^-$, $J^+$ are  linear operator acting on vector

\begin{eqnarray}\label{iter2}\nonumber
 J^-_j(u)(x)&=&i\int_{x_0}^x \sum_k \frac {V_{jk}(y)} {2k}  u_k(y) dy=-i\int 1_{y\geq x} \sum_k \frac {V_{jk}(y)} {2k}  u_k(y) dy\\
J^+_j(u)(x)&=&i\int_{x_0}^x e^{2ik(x-y)}
\sum_k \frac {V_{jk}(y)} {2k}    u_k(y)dy=-i\int 1_{y\geq x}  e^{2ik(x-y)}
\sum_k \frac {V_{jk}(y)} {2k}    u_k(y)dy.\end{eqnarray}

(This solve then the relation 
\begin{equation}\nonumber
\begin{array}{c}{d\over idx}  \gamma^-_{j,n+1} =
 \sum_k \frac {V_{jk}} {2k}   \gamma^-_{k,n}
+\sum_k -\frac {V_{jk}} {2k}  \gamma^+_{k,n}\\
{d\over idx}  \gamma^+_{j,n+1}-2k \gamma^+_{j,n+1} =
\sum_k \frac {V_{jk}} {2k}    \gamma^-_{k,n}
 +\sum_k -\frac {V_{jk}} {2k} \gamma^+_{k,n}.)
\end{array}
\end{equation}

We remark that the $\gamma^-_{j,n+1}(x)$ are constant in $x$ for $x\geq x_0$.
The equation could be see also as

$$ \gamma^-_{j,n+1}(x)- \gamma^+_{j,n+1}(x)=(J^-_j-J^+_j)(\gamma^-_{n}-\gamma^+_{n})(x),$$
$$ \gamma^+_{j,n+1}(x)=J^+_j(\gamma^-_{n}-\gamma^+_{n})(x).$$
So we get the expression
$$ \gamma^-_{j,n}(x)- \gamma^+_{j,n}(x)=[(J^--J^+)^n]_j (\gamma^-_{0}-\gamma^+_{0})(x), \mbox{ for } n\geq 0$$
$$ \gamma^+_{j,n}(x)=J^+_j (J^--J^+)^{n-1}
(\gamma^-_{0}-\gamma^+_{0})(x)\mbox{ for } n\geq 1,$$
$$ \gamma^-_{j,n}(x)=J^-_j (J^--J^+)^{n-1}  (\gamma^-_{0}-\gamma^+_{0})(x)\mbox{ for }
n\geq 1.$$

\begin{theorem}
Let $\im k\leq 0$, $|k|$ large enought,
the  formal series $\sum_n\gamma^\pm_{j,n}$ are absolutely convergent .
We have for each $n\in\N$ and $N\in \N$
\begin{equation}\label{ti1}
\sum_{n=0}^\infty \gamma^\pm_{j,n}-\sum_{n=0}^N \gamma^\pm_{j,n}=\O(|k|^{-N-1}),
\end{equation}
The series $ \sum_{n=N}^{\infty}(J^--J^+)^{n}$ are convergent   and theirs kernels, denote $k_N(x,y)$, satisfy for $N\geq 1$ $k_N(x,y)=B_N(k,x,y)V(y)$ with $B_N(k,x,y)$ is a matrix holomorphic in $k$ that satisfy

\begin{equation}\label{ti2}
|B_N(k,x,y)|\leq \frac {\|V\|_{L^1}^{N-1}} { {(2|k|)}^N }
 \end{equation}
umiformly in $x$.
\end{theorem}
Proof~:

Let $\im k\leq 0$.
We get $| e^{2ik(x-y)}|\leq 1$ for all $y>x$.

But the kernel of $J^--J^+$ is 
$$k(x,y)= -i 1_{y\geq x} (1- e^{2ik(x-y)} ) \frac {V(y)} {2k}$$
 which satisfy 
$|k(x,y)|\leq \frac {|V(y)|} {2k}$ uniformly in $x$.
So $\|J^--J^+\|_\infty\leq  \frac 1 {|2k|}   \int_{-\infty}^{\infty}{ \sum_{p,m} |V_{p,m}|(y)dy}.$

The kernel of $J^-$ is 
$$k^-(x,y)= -i 1_{y\geq x} \frac {V(y)} {2k}$$
 which satisfy 
$|k^-(x,y)|\leq \frac {|V(y)|} {2k}$ uniformly in $x$.
So $\|J^-\|_\infty\leq  \frac 1 {|2k|}   \int_{-\infty}^{\infty}{ \sum_{p,m} |V_{p,m}|(y)dy}.$

The kernel of $J^+$ is 
$$k^+(x,y)= -i 1_{y\geq x}  e^{2ik(x-y)} \frac {V(y)} {2k}$$
 which satisfy 
$|k^+(x,y)|\leq \frac {|V(y)|} {2k}$ uniformly in $x$.
So $\|J^+\|_\infty\leq  \frac 1 {|2k|}   \int_{-\infty}^{\infty}{ \sum_{p,m} |V_{p,m}|(y)dy}.$ 
Form which the results (\ref{ti1}) and (\ref{ti2}) holds.

\begin{remark}
The equation (\ref{eq:gen4}) is equivalent to

\begin{equation}\label{eqgen3b} 
 \left(\begin{array}{cc}
{d\over idx}&0\\
0&{d\over idx}-2k
\end{array}\right)  \left(\begin{array}{c} \gamma^--\gamma^+ \\ \gamma^+\end{array}\right)  \,
=\left(\begin{array}{cc} 0 & -2k  \\
 \frac V {2k} & 0   \end{array}\right)
\left(\begin{array}{c} \gamma^--\gamma^+ \\ \gamma^+\end{array}\right)
.\end{equation}

We note, that the preceding equations for $\gamma^\pm_{j,n}$ are similar
to the ones obtained by an exact WKB construction for scalar Schr\"odinger
equations, see for example the work of C.~Gerard and A.~Grigis~\cite{GG}
or T.~Ramond~\cite{R}.
\end{remark}

\section{the transmission matrix}\label{trans}

Let $\im k\leq 0$, $k\neq 0$, and $c={\,}^t(c_1,\cdots,c_n)$. 
Let  $v$ a solution of the equation (\ref{eq:gen1}) that behave like 
 $$v(x)=\sum_{j=1}^n \gamma^-_j(x,c) e^{-ikx} e_j+ \gamma^+_j(x,c)e^{-2ikx}  e^{ikx} e_{j+n},$$
with $\gamma^-_j(\infty,c)= c_j$ and $ \gamma^+_j(\infty,c)=0$  and
$\gamma^\pm_j(x,c)$ are solution of (\ref{eq:gen5}) (here $x_0=\infty$).\par
let $d={\,}^t(d_1,\cdots,d_n)$.
Let $w$ be solution  of the equation (\ref{eq:gen1}) that behave like
$$w(x)=\sum_{j=1}^n \gamma^-_j(x,de^{2iky}) e^{-ikx} e_j+
\gamma^+_j(x,de^{2iky})e^{-2ikx}  e^{ikx} e_{j+n},$$
with $\gamma^-_j(\infty,de^{2iky})= 0$ and $ \gamma^+_j(\infty,de^{2iky})=d e^{2iky}.$

The transmission matrix of the system (\ref{eq:gen1}) with is also the one of
the system (\ref{eq:gen0}) and (\ref{eq:gen3})  
is by definition the matrix of the linear map  associating to ${\,}^t(c,d)$ the vector
$$\begin{array}{ll} \lim\limits_{x\to-\infty}\gamma^-_j(x,c)+\lim\limits_{x\to-\infty} \gamma^-_j(x,de^{2iky})\\
\lim\limits_{x\to-\infty} \gamma^+_j(x,c)e^{-2ikx}+\lim\limits_{x\to-\infty}\gamma^+_j(x,de^{2iky})e^{-2ikx}.\end{array}$$
The transmission matrix is  analytic in  $k$.
Let denote the  transmission  matrix  
$$T(k)=\left(\begin{array}{ll} \tau_{11}(k)&\tau_{12}(k)\\ \tau_{21}(k)&\tau_{22}(k)\end{array}\right),$$
where $\tau_{ij}$ are $n\times n$ matrix given by the expression

$$\tau_{11}(k)(c)= c +\sum_{n=1}^\infty \lim_{x\to -\infty} J^-
(J^--J^+)^{n-1} (c)$$
$$\tau_{12}(k)(d)= -\sum_{n=1}^\infty \lim_{x\to -\infty} J^-(J^--J^+)^{n-1} (de^{2iky})$$
$$\tau_{21}(k)(c)= \sum_{n=1}^\infty \lim_{x\to -\infty} e^{-2ikx}J^+ (J^--J^+)^{n-1} (c)$$
$$\tau_{22}(k)(d)= d -\sum_{n=1}^\infty \lim_{x\to -\infty} e^{-2ikx}J^+ (J^--J^+)^{n-1} (de^{2iky})$$

For $k\in \R$, and $d\in\R^n$,
 $ \left(\begin{array}{ll}0&I_n\\I_n&0\end{array} \right) \overline{ w}$ is a solution of
 (\ref{eq:gen1}) (we use here $\overline V=V$) that behave like 
$$\overline w(x)=\sum_{j=1}^n \overline\gamma^-_j(x,de^{2iky})e^{2ikx} e^{-ikx} e_{j+n}+
\overline\gamma^+_j(x,de^{2iky})e^{2ikx} e^{-ikx} e_{j},$$
with $\overline\gamma^-_j(\infty,de^{2iky})= 0$ and $\overline
\gamma^+_j(\infty,de^{2iky})e^{2ikx}=d $ .
This proove that
$$\overline{\gamma^+_j(x,de^{2iky})e^{-2ikx}}=\gamma^-_j(x, d)$$
i.e. $\tau_{22}(k)= \overline{\tau_{11}( k)}$ for $k\in \R$, 
so $\tau_{22}(k)= \overline{\tau_{11}({\overline  k})}$ for $k\in\C$.
we get also using the same symetry $\tau_{12}(k)=
\overline{\tau_{21}(\overline k)}$.

\section{symmetry of the problem}\label{syme}

The scattering matrix is defined form the transmission matrix by the relation
$$S=\left(\begin{array}{ll} \tau_{11}- \tau_{12}\tau_{22}^{-1} \tau_{21} &  \tau_{12} \tau_{22}^{-1}\\
- \tau_{22}^{-1}\tau_{21}&\tau_{22}^{-1}\end{array}\right).$$
We remark that
$$S=\left(\begin{array}{ll} \tau_{12} & 0\\ 0& \tau_{22}^{-1}\end{array}\right)
\left(\begin{array}{ll}\tau_{12}^{-1}\tau_{11} &\tau_{22}^{-1} \\
    0&I  \end{array}\right)
\left(\begin{array}{ll} I & 0\\ -\tau_{21}& I \end{array}\right).$$
So 
\begin{equation}
\label{sca}
\det S(k)=\det \tau_{11} \times \det  \tau_{22}^{-1}=\frac {  {\det
    \tau_{11}( k)}}{\overline  \det \tau_{11}(\overline k)} .
\end{equation}
In particular 
\begin{equation}\label{dets}| \det S(k)|=1 \mbox{  for }k\in\R. \end{equation}
The zero of $\det S(k)$ are the complex conjugate of the pole of $\det S$. 
The pole of  $\det S$ are the resonances of our problem so we just want to count the number of zero of  $\det S(k)$ in the upper half plan. The function   $\det S(k)$ is holomorphic in the upper half plan so we going to use the theorem \ref{levin}
We need only to show  that $\det S$ is a function of exponential type in the full plan, the other hypothesis is satisfy using (\ref{dets}). 
By (\ref{sca}) the type of $\det S$ in a direction  is given by the type of $\det  \tau_{11}(k)$ minus the type of 
$\det \tau_{22}(k) $ in the same direction.   

Let define a Wronskian for our problem in some situation. Let $P$ be a projector on a space of dimension $p$ and  $\mathcal{ P}$ the projector defined by   $\mathcal{ P} \left(\begin{array}{l} u_1\\u_2\end{array}\right)= \left(\begin{array}{l} Pu_1\\Pu_2\end{array}\right)$. Denote by $P$ the matrix of $P$.
$W_P(u,v)(x)=\,^{t}(\mathcal{ P}u(x)) \left(\begin{array}{ll} 0&I\\-I&0\end{array}\right )\mathcal{ P}v(x)$
If $u$ and $v$ are solution of the system (\ref{eq:gen1}), and if $PV=\,^{t}VP$  then the Wronskian is constant in $x$.
This imply that if $PV=\,^{t}VP$ the transmission  matrix satisfy the relation
$$\,^tT(k)      
\left(\begin{array}{ll} 0&\,^{t}P P\\-\,^{t}P P&0\end{array}\right)  T(k)=
\left(\begin{array}{ll} 0&\,^{t}P P\\-\,^{t}P P&0\end{array}\right).$$
So we obtained $|\det    \mathcal{ P} T(k) \mathcal{P} |=1$.

\par\noindent\section{exponential type of the scattering matrix}\label{scatte}
\par\noindent
Using the  result (\ref{ti1}) for $\Im k<0$ and $|k|$ big enought, we obtain 
$$\tau_{11}(k)(c)= c+\O(\frac 1 {k}),$$ so $\det \tau_{11}(k)$ is of type $0$
in the lower half plan.
We have  
\begin{equation}\nonumber\begin{array}{lll}
\tau_{22}(k)(d)&=&d 
+\frac i {2k} \int V(z)  d dz \\
&+& \frac 1 {(2k)^2} \int  \int e^{-2ikz} V(z) 1_{y>z} V(y) d e^{2iky} dy dz\\
&-& \frac 1 {(2k)^2} \int  \int  V(z) 1_{y>z}  V(y) d dy dz\\
&+& \frac 1 {(2k)^2} \int\int  e^{-2ikz}  V(z)   B_2(z,y)  V(y)  de^{2iky} dy dz .\end{array}
\end{equation}

The type of $$d +\frac i {2k} \int V(z)  d dz- \frac 1 {(2k)^2} \int  \int  V(z) 1_{y>z}  V(y) d dy dz$$ is zero in the lower half plan.
The type of  $$ \int  \int e^{-2ikz} \sum_{p} V_{ip} (z) 1_{y>z}  V_{pj}(y) d e^{2iky} dy dz$$
is  in the lower half plan at most $\sup_{p} -u_{i,p}+t_{p,j}$ in the lower half plan.

The type of  $$\frac 1 {(2k)^2} \int\int  e^{-2ikz} \sum_{p,p'}  V_{i,p}(z)   B_{2,p,p'}(z,y)  V_{p',j}(y)  de^{2iky} dy dz $$ is at most $\sup_{p,p'} -u_{i,p}+t_{p',j}$ in the lower half plan.
Taking the determinant this give that 
$\det \tau_{22}(k)$ is of exponential type in the lower half plan with type at most 
$$\sup_{\sigma\in S_n}\sum_{i}  \sup_{p,p'} -u_{i,p}+t_{p',\sigma(i)}.$$

 We also get the relation ${\overline \det S(\overline k)}=\frac 1 {\det S(k)}$ form the expression (\ref{sca}).
So $\det S$ is of exponential type in $\{\im k>0\}$, in the real axis also by (\ref{dets}), as a conclusion $\det S$ is of exponential type in all the plan. So we have proof the first part of Theorem \ref{main}.
\par
Now let study the type of  $\det S$ in the direction $-i\R^+$,  this type is the constant $\pi C(V)$ of Theorem \ref{main}. 

\underline{Proof of theorem \ref{main} if $V$ is triangular}
Let write $V(x)=D(x)+N(x)$ then any product of the form $V(x_1) V(x_2)...V(x_n)$ is equal to 
$D(x_1)D(x_2)...D(x_n)$ plus a nilpotent matrix $N(x_1,\cdots,x_n)$. In particular this apply to the kernel of product of operators $J^+$ $J^-$ and $J^--J^+$. 
Comparing   the transition matrix for $V$ and $D$  this give that they exist nilpotent matrix $N_{11}, N_{12},N_{21},N_{22}$ such that 
$$\tau_{11}(V)(k)(c)=\tau_{11}(D)(k)(c)+N_{11}(k)$$
$$\tau_{22}(V)(k)(c)=\tau_{22}(D)(k)(c)+N_{22}(k)$$
$$\tau_{12}(V)(k)(c)=\tau_{12}(D)(k)(c)+N_{12}(k)$$
$$\tau_{21}(V)(k)(c)=\tau_{21}(D)(k)(c)+N_{21}(k)$$
So $$\det \tau_{11}(V)(k)=\det \tau_{11}(D)(k)(c)$$
and $$\det \tau_{11}(V)(k)=\det \tau_{11}(D)(k)(c).$$ 
So the type for the matrix $V$ is the same as the type for a diagonal matrix $D$. Then using the result of \cite{froese} or using  the wronskien we get 
$$  \tau_{22}(D)_{ii}\tau_{22}(D)_{ii}   - \tau_{12}(D)_{ii}\tau_{21}(D)_{ii}   =1.$$ 
The type of $(\tau_{22})_{ii}$  in the lower half plan is the same as the type of  $(\tau_{12})_{ii}(\tau_{21})_{ii} $. Using the expression of $\tau_{12}$, $\tau_{21}$ and using Lemma \ref{paley}  we find that the type  of $\tau_{22}(D))_{ii}$ is  $-u_{i,i}+t_{i,i}.$
So $$C(V)=C(D)= \frac 1 \pi \big(\sum_{i} u_{i,i}-t_{i,i}\big).$$

\underline{Proof of theorem \ref{main} if $V$ satisfy {\bf (H2)}.}
 We have  
\begin{equation}\nonumber\begin{array}{lll}
(\tau_{22}(k))_{ij}&=&\delta_{ij} 
+\frac i {2k} \int V_{ij}(z)  dz -\frac 1 {(2k)^2} \int  \int  \sum_p V_{i,p}(z) 1_{y>z}  V_{p,j}(y)   dy dz\\
&+& \frac 1 {(2k)^2} \int  \int e^{-2ikz} \sum_p V_{i,p}(z) 1_{y>z} (1+B_{2,p,p}(k,z,y)) V_{p,j} (y)  e^{2iky} dy dz\\
&+& \frac 1 {(2k)^2} \int\int  e^{-2ikz}  \sum_{p\neq p'} V_{i,p}(z)   B_{2,p,p'}(k,z,y)  V_{p',j} (y)  e^{2iky} dy dz .\end{array}
\end{equation}

The type of 
$$\frac 1 {(2k)^2} \int  \int e^{-2ikz}  V_{i,p(i)}(z) 1_{y>z} (1+B_{2,p,p}(k,z,y)) V_{p(i),\sigma_0(i)} (y)  e^{2iky} dy dz$$ is at least $t_{i,p(i)}-u_{p(i),\sigma_0(i)}$ in the direction $-i\R^+$. 
The proof is contain in the inequality, let $k=-ik'$ with $k'>0$,
\begin{equation}
\begin{array}{l} 
\int_z \int_y e^{-2k'z}  V_{i,p}(z) 1_{y>z} (1+B_{2,p,p}(k,z,y)) V_{p,j} (y)  e^{2iky} dy dz\\
\geq  
\int _{z\leq  t_{i,p}+\epsilon} \int_{y\geq  u_{p(i),\sigma_0(i)}-\epsilon} \ e^{-2k' ( t_{i,p}+\epsilon) }  V_{i,p}(z) \frac 1 2  V_{p,j} (y)  e^{2k' (u_{p(i),\sigma_0(i)}-\epsilon)} dy dz\\
\geq  e^{-2k' ( t_{i,p}+\epsilon) }  e^{2k' (u_{p(i),\sigma_0(i)}-\epsilon)} C(\epsilon)
\end{array}
\end{equation}
with $C(\epsilon)>0$ and independent of $k$.
The type of  
$$\frac 1 {(2k)^2} \int\int  e^{-2ikz} \sum_{p\neq p'}  V_{i,p}(z)   B_{2,p,p'}(k,z,y)  V_{p',j}(y)  de^{2iky} dy dz $$ is at most $\sup_{p\neq p'}t_{i,p}-u_{p',j}$ in the lower half plan.
Using the hypotheses {\bf(H2)} this give that 
$\det \tau_{22}(k)$ is of exponential type in the lower half plan with type  
$$\sum_{i}  \sup_{p} t_{i,p}-u_{p,\sigma_0(i)}.$$

\bibliographystyle{alpha}
\bibliography{abib}
\end{document}